\documentclass[11pt]{amsart}

\usepackage{epsf,amssymb,times,overpic,amscd}
\usepackage[usenames]{color}
\usepackage[all]{xy}

\newtheorem{thm}{Theorem}[section]
\newtheorem{prop}[thm]{Proposition}
\newtheorem{lemma}[thm]{Lemma}

\theoremstyle{definition}

\theoremstyle{remark}

\newtheorem{rmk}[thm]{Remark}

\newcommand{\R}{\mathbb{R}}
\newcommand{\Z}{\mathbb{Z}}

\newcommand{\N}{\mathbb{N}}

\newcommand{\bdry}{\partial}
\newcommand{\s}{\vskip.1in}
\newcommand{\n}{\noindent}

\newcommand{\F}{\mathbb{F}}

\newcommand{\op}{\operatorname}

\begin{document}

\title{$HF=ECH$ via open book decompositions: a summary}

\author{Vincent Colin}
\address{Universit\'e de Nantes, 44322 Nantes, France}
\email{Vincent.Colin@univ-nantes.fr}

\author{Paolo Ghiggini}
\address{Universit\'e de Nantes, 44322 Nantes, France}
\email{paolo.ghiggini@univ-nantes.fr}
\urladdr{http://www.math.sciences.univ-nantes.fr/\char126 Ghiggini}

\author{Ko Honda}
\address{University of Southern California, Los Angeles, CA 90089}
\email{khonda@usc.edu} \urladdr{http://math.usc.edu/\char126 khonda}

\date{This version: March 6, 2011.}

\keywords{contact structure, Reeb dynamics, embedded contact homology, open book decompositions, Heegaard Floer homology}

\subjclass[2000]{Primary 57M50; Secondary 53D10,53D40.}


\begin{abstract}
We sketch the proof of the equivalence between the hat versions of Heegaard Floer
homology and embedded contact homology, abbreviated ECH. The key point is to
express these two Floer homology theories in terms of an open book decomposition
of the ambient manifold.
\end{abstract}

\maketitle

\section{Introduction}

There are three Floer homology theories of a closed oriented $3$-manifold --- Heegaard Floer homology, embedded contact homology and Seiberg-Witten Floer homology --- which are conjectured to be equivalent.\footnote{There is a fourth one --- instanton Floer homology --- whose relationship to the other three Floer theories is not well understood.} Heegaard Floer homology was defined by Ozsv\'ath-Szab\'o~\cite{OSz1,OSz2}, embedded
contact homology (abbreviated ECH) by Hutchings~\cite{Hu1} and Hutchings-Taubes~\cite{HT1,HT2}, and Seiberg-Witten Floer homology by Kronheimer-Mrowka~\cite{KM2}.  These theories have had spectacular applications over the last decade, ranging from the Gordon conjecture due to Kronheimer-Mrowka-Ozsv\'ath-Szab\'o~\cite{KMOS} and progress in the exceptional surgery problem due to Ghiggini~\cite{ghiggini:fibered} and Ni~\cite{ni:fibered} to the classification of Seifert manifolds admitting tight contact structures due to Lisca and Stipsicz~\cite{lisca-stipsicz} and the Weinstein conjecture in dimension three due to Taubes~\cite{T1}, just to name a few.

In 2006, Taubes proved a breakthrough compactness theorem in Seiberg-Witten theory, enabling him to prove the Weinstein conjecture in dimension three~\cite{T1} and to establish the equivalence of Seiberg-Witten Floer cohomology and ECH~\cite{T2}, shortly thereafter. In order to establish the equivalence of all three theories, it then suffices to show that Heegaard Floer cohomology is equivalent to either ECH or Seiberg-Witten Floer cohomology.  Since Heegaard Floer cohomology and ECH are both Floer theories which involve holomorphic curves, it is natural to expect a chain map between the two; this is indeed our approach.

The goal of this research announcement is to present some ingredients of the proof of the following theorem:

\begin{thm}\label{thm: isomorphism}
Let $M$ be a closed oriented $3$-manifold, $\xi$ a positive cooriented contact structure on $M$ and $\mathfrak{s}_\xi$
the associated $\mbox{Spin}^c$-structure. Then $$\widehat{HF}(-M,\mathfrak{s}_\xi+PD(A))\simeq \widehat{ECH}(M,\xi,A),$$ where $A\in H_1(M;\Z)$.
\end{thm}

Here $\widehat{HF}$ and $\widehat{ECH}$ are the hat versions of Heegaard Floer homology and ECH, respectively; their definitions will be reviewed shortly. Also $-M$ is $M$ with the opposite orientation, and $\widehat{HF}(-M)$ is isomorphic to the Heegaard Floer hat cohomology group of $M$. For simplicity we assume that our coefficient system is $\F=\Z/2\Z$.  The full details of Theorem~\ref{thm: isomorphism} will appear in \cite{CGH1,CGH2}.

Using the same range of ideas and a spectral sequence argument, we expect to reduce the plus case to the hat case and prove the following isomorphism:
$$HF^+ (-M ,\mathfrak{s}_\xi +PD(A))\simeq ECH (M,\xi ,A).$$
This is work in progress~\cite{CGH3}.

\begin{rmk}
A proof of the isomorphism  $HF^+ (-M ,\mathfrak{s}_\xi +PD(A))\simeq ECH (M,\xi ,A)$
has been announced by Kutluhan-Lee-Taubes~\cite{KLT}, using different methods. From
their proof an alternative proof of the isomorphism between the hat versions follows by
some simple homological algebra.
\end{rmk}

\section{Review of Heegaard Floer homology}
\label{section: review of HF}

In this section we briefly review the Heegaard Floer homology groups associated to a closed oriented $3$-manifold $M$.

Heegaard Floer homology was introduced by Ozsv\'ath and Szab\'o~\cite{OSz1,OSz2} in an attempt to obtain a more combinatorial (and easier-to-define) version of the Seiberg-Witten Floer homology of $M$.\footnote{The program of obtaining a combinatorial version of Heegaard Floer homology was recently brought to fruition by Sarkar-Wang~\cite{SW} for the hat version and by Manolescu-Ozsv\'ath-Thurston~\cite{MOT} for the plus version.} The definition of $\widehat{HF}(M)$ that we give in this paper is not the original definition in \cite{OSz1,OSz2}, but rather a reformulation proposed by Eliashberg and carried out by Lipshitz~\cite{Li}.

The group $\widehat{HF}(M)$ is defined via a {\em pointed Heegaard diagram} $(\Sigma,\alpha,\beta,z)$ associated to a Heegaard decomposition of $M$. Here $\Sigma$ is a closed, oriented, connected surface of genus $g$ which divides $M$ into two handlebodies $H_\alpha$ and $H_\beta$, $\alpha =\{\alpha_1,\dots,\alpha_{g}\}$ (resp.\  $\beta =\{\beta_1,\dots,\beta_{g}\}$) is a collection of $g$ pairwise disjoint simple closed curves which bound disks in $H_\alpha$ (resp.\ $H_\beta$) and are linearly independent in $H_1(\Sigma;\Z)$, and $z\in\Sigma-\alpha-\beta$.

Given $y \in \Sigma$ we will use the symbol $\vec{y}$ to denote the chord $[0,1] \times \{ y \} \subset [0,1] \times \Sigma$. The chain group $\widehat{CF} (\Sigma,\alpha,\beta,z)$ is the free $\F$-vector space
generated by $\mathcal{S}_{\alpha,\beta}$, the set of $g$-tuples of chords $\{ \vec{y}_1,
\dots, \vec{y}_g \}$ in $[0,1] \times \Sigma$ from $\{0\}\times \beta$ to $\{1\}\times
\alpha$, for which there exists a permutation $\sigma\in \mathfrak{S}_{g}$ satisfying
$y_i \in \alpha_i\cap \beta_{\sigma (i)}$, $i=1,\dots,g$. We will write $\vec{\bf y}$ for the
$g$-tuple of chords $\{\vec{y}_1, \ldots , \vec{y}_g \}$.
In order to define the differential $\bdry \vec{\bf y} =
\sum_{\vec{\bf y}'\in \mathcal{S}_{\alpha,\beta}} \langle \bdry \vec{\bf y}, \vec{\bf y}'
\rangle \vec{\bf y}'$, we consider the symplectic fibration
$$\pi: (W= \R \times [0,1] \times \Sigma,\Omega) \to (\R\times[0,1],ds\wedge dt),$$
where $s,t$ are the coordinates on $\R \times [0,1] $, $\omega$ is an area form on
$\Sigma$, and $\Omega =ds \wedge dt +\omega$ on $W$.  The submanifolds
$C_\alpha =\R \times \{ 1\}\times \alpha$ and $C_\beta =\R \times \{0\} \times \beta$
are Lagrangian submanifolds of $(W,\Omega)$ which are contained in $\partial W$.
Next, let $J$ be an almost complex structure on $W$ which is tamed by $\Omega$, is
$s$-invariant, and sends $\bdry_s\mapsto \bdry_t$ and $T\Sigma$ to itself. Then
$\langle \bdry \vec{\bf y}, \vec{\bf y}'\rangle$ is a mod $2$ count of embedded, degree
$g$, Fredholm index $\op{ind}=1$, $J$-holomorphic multisections of $\pi$ which are
asymptotic to the chords $\vec{\bf y}$ as $s\to +\infty$ and to the chords $\vec{\bf y}'$
as $s\to -\infty$, and whose boundary maps to $C_\alpha$ and $C_\beta$.\footnote{In \cite{Li}, Lipshitz uses $W= \Sigma \times [0,1] \times \R$, while we use $W=  \R \times [0,1] \times \Sigma$. This accounts for some slight differences, e.g., in \cite{Li}, $\langle \bdry \vec{\bf y}, \vec{\bf y}'\rangle$ is a count of curves which are asymptotic to $\vec{\bf y}$ as $s\to -\infty$ and to $\vec{\bf y}'$ as $s\to +\infty$, whereas our asymptotics are reversed.} The homology of
the chain complex $\widehat{CF}(\Sigma,\alpha,\beta,z)$ is independent of the choice of pointed Heegaard diagram for $M$, and will be written as $\widehat{HF}(M)$.

The Heegaard Floer groups can be decomposed along Spin$^c$-structures on $M$:
$$\widehat{HF} (M)=\bigoplus_{\mathfrak{s} \in \op{Spin}^c (M)} \widehat{HF} (M,\mathfrak{s} ).$$
See~\cite[Section~2]{Li} for more details.

\begin{rmk}
A {\em stable Hamiltonian structure} $(\lambda,\omega,R)$ on an oriented $3$-manifold consists of a $1$-form $\lambda$, a closed nowhere zero $2$-form $\omega$, and a vector field $R$ (called the {\em Hamiltonian vector field}) such that  $\lambda(R)=1$, $R$ directs $\ker \omega$, and $\ker \omega\subset \ker d\lambda$.  On $[0,1]\times \Sigma$ with coordinates $(t,x)$, consider the stable Hamiltonian structure $(dt, \omega, \bdry_t)$, where $\omega$ is (the pullback of) an area form on $\Sigma$.  Then the chords $\vec{\bf y}$ are $g$-tuples of chords of the Hamiltonian vector field $\bdry_t$. In this way, $\widehat{CF}(\Sigma,\alpha,\beta,z)$ can be placed in the context of symplectic field theory~\cite{BEHWZ,EGH}.
\end{rmk}

\section{Review of embedded contact homology}
\label{section: review of ECH}

In this section we briefly review the embedded contact homology (ECH) groups associated to a closed oriented $3$-manifold $M$.  ECH was defined by Hutchings~\cite{Hu1} and Hutchings-Taubes~\cite{HT1,HT2} and is intimately connected with the dynamics of a Reeb vector field.

A {\em (positive) contact form} $\lambda$ on $M$ is a $1$-form satisfying $\lambda \wedge d\lambda >0$. The {\em Reeb vector field} $R=R_\lambda$ of $\lambda$ is given by $i_R d\lambda=0$ and $i_R\lambda=1$.  We assume that $R$ is {\em nondegenerate}, i.e., the first return map along each (not necessarily simple) periodic orbit does not have $1$ as an eigenvalue.

\subsection{$ECH(M)$}

The chain complex $ECC(M,\lambda,J)$ is the free $\F$-vector space generated by multisets (i.e., sets where elements are allowed to have multiplicities $\in\N$) of simple periodic orbits of $R_\lambda$, where the multiplicity assigned to a simple hyperbolic orbit is always $1$. Such multisets are called {\em orbit sets}.  Let $J$ be an almost complex structure which is adapted to the symplectization $(\R\times M,d(e^s\lambda))$, i.e., sends $\xi=\ker \lambda$ to itself and $\bdry_s$ to $R_\lambda$, where $\bdry_s$ is the $\R$-coordinate. Let $\gamma$ and $\gamma'$ be orbit sets of $R_\lambda$. Then the ECH differential $\langle \bdry \gamma,\gamma'\rangle$ is the (mod 2) count of ECH index $I_{ECH}=1$ holomorphic curves $u$ in $\R\times M$ {\em from $\gamma$ to $\gamma'$}. (By this we mean that $u$ is asymptotic to $\gamma$ as $s\to +\infty$ and to $\gamma'$ as $s\to -\infty$, with the correct multiplicity.)  The essential ingredient in the definition is the {\em ECH index} $I_{ECH}$ of a relative homology class $Z\in H_2(M,\gamma\cup\gamma')$.  Although we do not give its precise definition here, $I_{ECH}$ is roughly the sum of two terms: the Fredholm index and the self-intersection number of $Z$, including the asymptotic self-intersection as $s\to \pm \infty$. By the positivity of intersections of $J$-holomorphic curves in dimension four, $I_{ECH}=1$ curves are (for the most part) embedded.

The homology group $ECH(M,\lambda,J)$ turns out to be independent of the choices of $\lambda$ and $J$, and in particular is independent of the choice of contact structure on $M$.  There is still no direct proof of this invariance, and the only known proof is through Taubes' isomorphism~\cite{T2} between $ECH(M,\lambda,J)$ and Seiberg-Witten Floer cohomology.

The contact structure $\xi =\ker \lambda$ determines a unique Spin$^c$-structure $\mathfrak{s}_\xi$ and $ECH(M,\xi)$ can be decomposed using the total homology class $A\in H_1 (M;\Z)$ of $\gamma$, where the orbits are counted with the appropriate multiplicities.

\subsection{$\widehat{ECH}(M)$}

We now define the variant $\widehat{ECH}(M)$ of $ECH(M)$, called the {\em ECH hat group}. First pick a generic point $z\in \R\times M$ and take $\gamma$ and $\gamma'$ to be orbit sets of the $ECH$ chain complex. We define the map:
$$U: ECC(M,\lambda,J)\to ECC(M,\lambda,J),$$
where $\langle U(\gamma ),\gamma'\rangle$ is the (mod 2) count of holomorphic curves of ECH index $I=2$ from $\gamma$ to $\gamma'$ which pass through the point $z$. The ECH hat group $\widehat{ECH}(M,\lambda,J)$ is then defined as the mapping cone of $U$. The group $\widehat{ECH}(M,\lambda,J)$ also has an interpretation as a sutured ECH group by the work of Colin-Ghiggini-Honda-Hutchings~\cite{CGHH}.

\section{Open book decompositions}

Let $M$ be a closed oriented $3$-manifold. An {\it open book decomposition} of $M$ is a triple $(S,h,\phi)$, where $S$ is a compact, oriented, connected surface with nonempty boundary, $h:S\stackrel\sim\rightarrow S$ is an orientation-preserving diffeomorphism such that $h|_{\bdry S}=id$, and $\phi$ is an orientation-preserving homeomorphism from $(S\times [0,1])/\sim$ to $M$. Here the equivalence relation $\sim$ is generated by $(x,1)\sim (h(x),0)$ for all $x\in S$ and $(y,t)\sim (y,t')$ for all $y\in \partial S$ and $t,t'\in [0,1]$. We refer to $(\partial S\times [0,1])/\sim$ as the {\it binding} and $S_t=S\times \{ t\}$ as a {\it page} of the open book decomposition.  The homeomorphism $\phi$ will usually be suppressed from the notation.

Let $(S,h)$ be an open book decomposition of $M$. Then a contact structure $\xi$ is {\em adapted to} $(S,h)$ if it admits a Reeb vector field which is positively transverse to the pages and is tangent to and directed by the binding (here the binding is oriented as the boundary of a page); such a contact structure will be denoted by $\xi_{(S,h)}$.  In the 1970's, Thurston and Winkelnkemper~\cite{TW} discovered a construction which assigns a contact structure $\xi_{(S,h)}$ to any open book decomposition $(S,h)$. More recently, in the fundamental work \cite{Gi2}, Giroux showed that this assignment gives rise to the following one-to-one correspondence:

\begin{thm}[Giroux, \cite{Gi2}]
There is a one-to-one correspondence between isotopy classes of contact structures on $M$ and isotopy classes of open book decompositions $(S,h)$ modulo positive stabilization.
\end{thm}

There is a related manifold, the {\it suspension} $N_{(S,h)}$ of $(S,h)$, which is obtained from $S\times [0,1]$ with coordinates $(x,t)$ by identifying $(x,1)\sim (h(x),0)$; we abbreviate it as $N$ when $(S,h)$ is understood.  Let $\pi: N\to S^1=[0,1]/\sim$ be the corresponding fibration. Since $h\vert_{\partial S} = id$, there is a natural oriented identification $\bdry N\simeq \R^2/\Z^2$ which sends $\partial S$, with the boundary orientation of $S$, to a closed curve directed by $(1,0)$ and $\{x\}\times[0,1]/\sim$, where $x\in \bdry S$ and the orientation is given by the usual orientation of $[0,1]$, to a closed curve directed by $(0,1)$. This identification allows us to refer to slopes of simple closed curves on $\bdry N$ or on tori parallel to $\bdry N$.

\section{How to define $\widehat{HF} (-M)$ from a page of an open book}
\label{section: defining HF from one page}

In this section we explain how to rephrase $\widehat{HF}(-M)$ in terms of an open book decomposition $(S,h)$ of $M$, using a construction discovered by \cite{HKM}. From now on we assume that $\bdry S$ is connected and $S$ has genus $g$.

The open book decomposition $(S,h)$ gives rise to a Heegaard decomposition $M=H_\alpha\cup H_\beta$, where $H_\alpha =S\times [0,{1\over 2}]/ \sim$, $H_\beta =S\times [{1\over 2}, 1]/ \sim$, and the Heegaard surface $\Sigma =S_{1/2} \cup -S_0$ is the union of two pages glued along the binding.

\begin{figure}[ht]
\begin{overpic}[width=6cm]{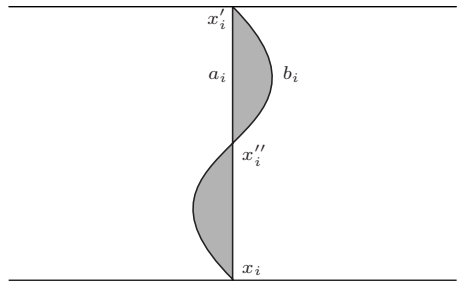}
\put(44.5,45){\tiny $a_i$} \put(61,45) {\tiny $b_i$} \put(52,27)
{\tiny $x_i''$} \put(52,2.1){\tiny $x_i$} \put(44.3,57.3){\tiny
$x_i'$}
\end{overpic}
\s
\caption{A portion of $S_{1/2}$. The shaded regions are the ``thin strips'' $D_i$ and
$D_i'$ which connect $x_i''$ to $x_i$ or $x_i'$.} \label{darkside}
\end{figure}

A {\it basis of arcs} for $S$ is a collection of $2g$ pairwise disjoint properly embedded arcs ${\bf a} =\{a_1,\dots,a_{2g}\}$ in $S$ such that $S-{\bf a}$ is a connected polygon. Starting from a basis ${\bf a}$ for $S$, we can construct $\alpha$- and $\beta$-curves for $\Sigma$ as follows: $\alpha = ({\bf a} \times \{{1\over 2}\}) \cup ({\bf a} \times \{ 0\})$ and $\beta = ({\bf b} \times \{{1\over 2}\}) \cup (h({\bf a} )\times \{ 0\})$. Here ${\bf b}$ is a small deformation of ${\bf a}$ relative to its endpoints, so that each pair $a_i$ and $b_i$ intersects each other transversely at three points: two of the intersections are their endpoints $x_i$ and $x_i'$ on $\bdry S$ and the third intersection is an interior point $x_i''$; see Figure~\ref{darkside}.
This means that all the intersection points of $\alpha$ and $\beta$ lie in $S_0$, with the
exception of the points $x_i''=x_i''\times\{{1\over 2}\}$.  We place the basepoint $z$ on
$S_{1/2}$, away from the ``thin strips'' $D_i$ and $D_i'$, $i=1,\dots,2g$, given in
Figure~\ref{darkside}. The positioning of $z$ prevents holomorphic curves involved in
the differential for $\widehat{CF} (-\Sigma,\alpha,\beta,z)$ --- besides ones corresponding to the thin strips --- from entering the $S_{1/2}$ region.  Hence all of the nontrivial holomorphic curve information is concentrated on $S_0$. Here $\widehat{HF}(-\Sigma, \alpha,\beta,z)$ is isomorphic to $\widehat{HF}(-M)$ since we reversed the orientation of $\Sigma$.

Let $\mathcal{S}_{{\bf a},h({\bf a})}\subset \mathcal{S}_{\alpha,\beta}$ consist of
$2g$-tuples $\vec{\bf y}$, all of whose intersections are in ${\bf a}\cap h({\bf a})$.
We then define $\widehat{CF}(S, {\bf a}, h({\bf a}))$ as the chain complex generated
by $\mathcal{S}_{{\bf a},h({\bf a})}$, modulo the identification $\{\vec{x}_i\}\cup\vec{\bf y}_0
\sim\{\vec{x}_i'\}\cup \vec{\bf y}_0$ for all $(2g-1)$-tuples of chords $\vec{\bf y}_0$,
and whose differential counts holomorphic curves in $\R\times [0,1]\times S$.  The following proposition was proved in \cite{CGH2}:

\begin{prop}\label{prop:hf}
$\widehat{HF}(S,\mathbf{a},h(\mathbf{a}),z)\simeq \widehat{HF}(-\Sigma,\alpha,\beta,z)$.
\end{prop}

\begin{rmk}
$\vec{\bf x} =\{\vec{x}_1,\dots,\vec{x}_{2g}\}$ is a cycle and its class $[\vec{\bf x}]
\in \widehat{HF} (S,\mathbf{a},h(\mathbf{a}))$ is the {\em contact class}
$c(\xi_{(S,h)})$ of $\xi_{(S,h)}$; see~\cite{HKM}.
\end{rmk}

\section{How to define $\widehat{ECH} (M,\xi)$ from a page of an open book}
\label{section: defining ECH from page}

Let $\xi=\xi_{(S,h)}$ be a contact structure which is adapted to the open book decomposition $(S,h)$ of $M$.  The goal of this section is to formulate $\widehat{ECH} (M,\xi)$ in terms of the suspension $\pi: N=N_{(S,h)}\to S^1$, i.e., to eliminate the binding. (See~Theorem \ref{thm: binding}.)  This is motivated in part by the work of Wendl~\cite{We} and Yau~\cite{Y}.

\subsection{A nice contact form}

The following lemma, proved in \cite{CGH1}, furnishes a nice contact form on $N$:

\begin{lemma}
There exists a contact form $\lambda'$ on $N$, whose Reeb vector field $R_{\lambda'}$ satisfies the following properties:
\begin{itemize}
\item $R_{\lambda'}$ is transverse to the fibers of $\pi$;
\item $R_{\lambda'}$ is tangent to $\partial_t$ along $\partial N$;
\item the first return map of $R_{\lambda'}$ on a small neighborhood $(\R/\Z) \times [-\varepsilon,0]$ of $\partial S=(\R /\Z)\times \{0\}$ in $S$ is $(x,r)\mapsto (x-r,r)$.
\end{itemize}
\end{lemma}

After a small perturbation of $\lambda'$, we obtain a contact form $\lambda$ whose Reeb vector field $R_\lambda$ is nondegenerate on $int(N)$ and is foliated by an $S^1$-Morse-Bott family of orbits along $\partial N$.  Pick two distinct orbits in the Morse-Bott family that we label as $e$ and $h$ and view as an elliptic and a hyperbolic orbit, respectively; see \cite{Bo1,Bo2} for more details on Morse-Bott theory.

\subsection{ECH groups on $N$}

Let $\mathcal{P}$ be the set of simple orbits of the Reeb vector field $R_\lambda$ in $int(N)$, together with orbits $e$ and $h$. Let $ECC(N,\lambda)$ be the ECH chain complex generated over $\F$ by orbit sets whose constituent simple orbits are in $\mathcal{P}$. The direct summand $ECC_j(N,\lambda)$ of $ECC(N,\lambda)$ consists of orbit sets whose total homology class has algebraic intersection number $j$ with a fiber $S$ of $\pi$. We write $ECH_j(N,\lambda)$ and $ECH(N,\lambda)$ for the homology of $ECC_j(N,\lambda)$ and $ECC(N,\lambda)$.

There are inclusions of chain complexes:
$$\mathfrak{I}_j: ECC_j(N,\lambda)\to ECC_{j+1}(N,\lambda),$$
given by $\gamma\mapsto e\gamma$, where the orbit set $\gamma$ is written multiplicatively as $\prod_i \gamma_i^{m_i}$. The collection of maps $\{\mathfrak{I}_j\}$ gives rise to the direct limit $\lim \limits_{\longrightarrow \atop i} ECH_i(N,\lambda)$.  We can view this direct limit as the homology of the quotient of $ECC(N,\lambda)$, obtained by identifying $e\gamma\sim\gamma$; the reason for this will be explained below.

\begin{rmk}
Given a contact structure $\xi$ on $M$, its contact class in $ECH(M,\xi)$ is represented by the empty set, written multiplicatively as $1$.
\end{rmk}

\subsection{Eliminating the binding}

The following was proved in \cite{CGH1}:

\begin{thm}\label{thm: binding}
$\displaystyle\widehat{ECH} (M,\xi) \simeq \lim \limits_{\longrightarrow \atop i} ECH_i(N,\lambda).$
\end{thm}

The intuitive idea of the proof is as follows. Let $\lambda_\delta$, $0 <\delta \ll 1$, be a contact form which is adapted to $(S,h)$ such that, on a solid torus neighborhood $D^2 \times (\R/\Z)$ of the binding $K$ with cylindrical coordinates $(r,\theta,z)$, its Reeb vector field $R_{\lambda_\delta}$
is tangent to the concentric tori $\{r=const\}$ and has constant slope $\frac{1}{\delta}$
away from $K$.  As we send $\delta\to 0$, the Conley-Zehnder index of the binding
goes to $+\infty$; we should therefore be able to ignore the binding as $\delta\to 0$.
At the same time, one sees that $I_{ECH}=1$ holomorphic curves that cross the binding
when $\delta>0$ are in one-to-one correspondence with holomorphic curves which
have $e$ at the negative end when $\delta=0$.  This is the reason for identifying $e=1$. Similarly, if we place the marked point $z$ on the binding, then $I_{ECH}=2$ holomorphic curves that pass through $z$ are in one-to-one correspondence with $I_{ECH}=1$ holomorphic curves which have $h$ at a negative end when $\delta=0$.

\subsection{Periodic Floer homology}

At this point, it is convenient to switch from the ECH groups to the similarly defined {\em periodic Floer homology} (PFH) groups, also defined by Hutchings~\cite{Hu1}. The relevant stable Hamiltonian structure $(\lambda,\omega,R)$ is given as follows: On $S\times[0,2]$ with coordinates $(x,t)$, let $\lambda=dt$, $R=\bdry_t$ and $\omega$ be an area form on $S$. We then identify $h: S\times\{2\}\stackrel\sim\to S\times\{0\}$, where $h^*\omega=\omega$, $h|_{\bdry S}=id$ and a certain technical condition called the {\em zero flux condition} is satisfied (see for example~\cite{CHL}).  The following is proved in \cite{CGH2}:

\begin{prop}
$ECH_i(N)\simeq PFH_i(N)$.
\end{prop}

\section{The map $\Phi:\widehat{CF}(S,\mathbf{a},h(\mathbf{a}))\to PFC_{2g}(N)$}

We define a map $\Phi : \widehat{CF}(S,\mathbf{a},h(\mathbf{a})) \to PFC_{2g} (N)$ by counting degree $2g$ embedded multisections in a symplectic fibration $(W_+,\Omega_+)$ which is a cylinder over $[0,1] \times S$ at the positive end and a cylinder over $N$ at the negative end. (At both ends we have the stable Hamiltonian structure $(dt,\omega,\bdry_t)$, where $\omega$ is an area form on $S$.)  This can be viewed as an amalgamation of the work of Seidel~\cite{Se} and Donaldson-Smith~\cite{DS}.

\subsection{The symplectic cobordism $(W_+,\Omega_+)$}

\begin{figure}[ht]
\begin{overpic}[width=8cm]{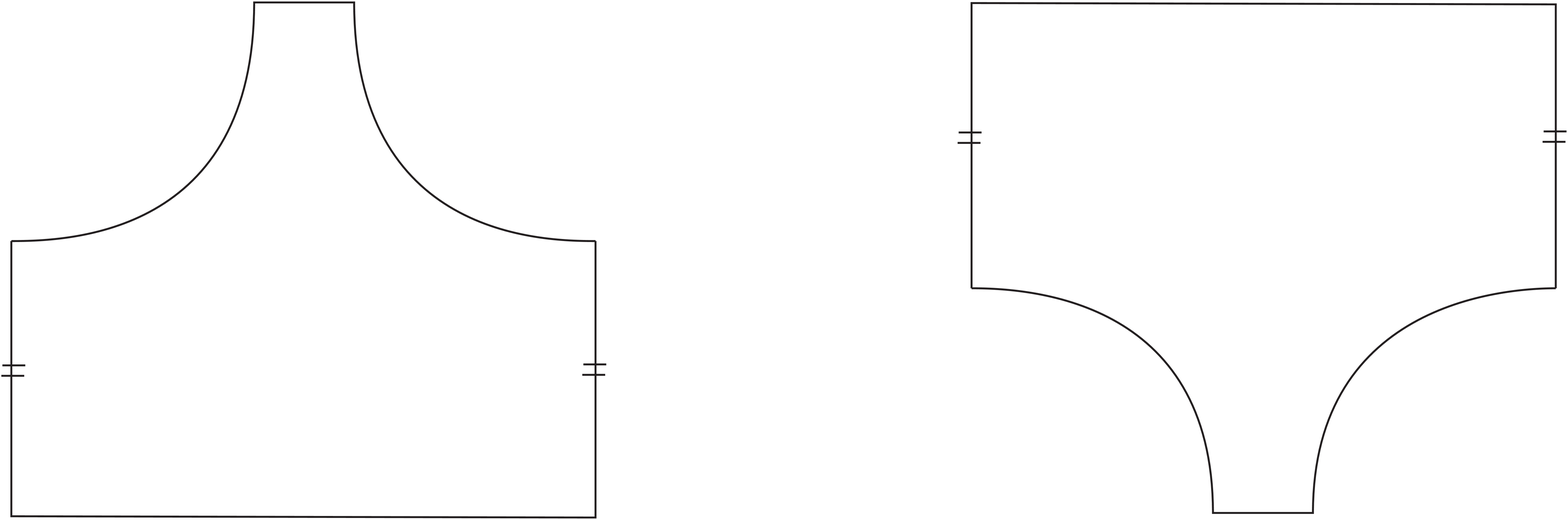}
\put(60.12,28){$\times$} \put(97.35,28){$\times$}
\end{overpic}
\s
\caption{The bases $B_+$ and $B_-$.  The sides are identified. Both $B_+$ and $B_-$ are biholomorphic to a disk with an interior puncture and a boundary puncture.}
\label{figure: bases}
\end{figure}

Consider the infinite cylinder $\R \times S^1= \R\times (\R/2\Z)$ with coordinates $(s,t)$. Let $\pi: N \to S^1$ be the fibration $(x,t)\mapsto t$ and let $\pi: \R \times N \to \R\times S^1$ be the extension $(s,x,t)\mapsto (s,\pi(x,t))$. We define $W_+ =\pi^{-1} (B_+)$, where $B_+^c$ is the subset $[2,\infty)\times[1,2]\subset \R\times(\R/2\Z)$ with the corners rounded and $B_+=(\R\times (\R/2\Z)) -B_+^c$. See the left-hand side of Figure~\ref{figure: bases}.
We then have a symplectic fibration $$\pi : (W_+ , \Omega_+=ds\wedge dt+\omega)\to (B_+,ds\wedge dt).$$ The fibration admits a symplectic connection defined as the $\Omega_+$-orthogonal to the tangent spaces of the fibers; it is spanned by $\partial_s$ and $\partial_t$ on $\R\times S\times[0,2]$ (before identifying $\R\times S\times\{2\}$ and $\R\times S\times\{0\}$).

\subsection{Lagrangian boundary}

Pick the point $(3,1) \in \partial B_+$ and consider the basis ${\bf a}$ in the fiber $\pi^{-1} (3,1)$. The Lagrangian submanifold $L_{\bf a}\subset (W_+ ,\Omega_+ )$ is defined as the trace of the parallel transport of ${\bf a}\subset \pi^{-1}(3,1)$ along $\partial B_+$ using the symplectic connection. Since the symplectic connection is spanned by $\partial_s$ and $\partial_t$ on $\R\times S\times[0,2]$, over the strip $\{ s\geq 3, t\in [0,1 ] \}$ we have:
$$L_{\bf a} \cap \{s\geq 3, t=0\} = \{ s\geq 3\} \times h({\bf a}) \times \{0\},$$
$$L_{\bf a} \cap \{s\geq 3, t=1\} =\{ s\geq 3\} \times {\bf a} \times \{ 1\}.$$

\subsection{The map $\Phi$}

Let $J_+$ be an almost complex structure on $(W_+,\Omega_+)$ which is the restriction of an adapted almost complex structure on $\R\times N$, i.e., takes $\bdry_s\mapsto \bdry_t$ and $TS$ to itself. Let $\vec{\bf y}$ be a $2g$-tuple of chords in $\widehat{CF}(S,{\bf a}, h({\bf a}))$ and $\gamma$ be an orbit set in $PFC_{2g}(N)$. If we write $\Phi(\vec{\bf y})= \sum_{\gamma}  \langle \Phi ( \vec{\mathbf{y}} ),\gamma \rangle \gamma$, then $\langle \Phi ( \vec{\mathbf{y}} ),\gamma \rangle$ counts degree $2g$, ECH index $I_{W_+}=0$ (briefly explained below), $J_+$-holomorphic multisections of the fibration $\pi: W_+\to B_+$ which are positively asymptotic to $\vec{\bf y}$ and negatively asymptotic to $\gamma$ and whose boundary is contained in distinct components of the Lagrangian boundary $L_{\bf a}$.

\begin{thm} \label{thm: Phi chain map}
The map $\Phi : \widehat{CF}(S,\mathbf{a},h(\mathbf{a})) \to PFC_{2g} (N)$ is a chain map.
\end{thm}

Theorem~\ref{thm: Phi chain map} is proved using standard arguments in symplectic geometry such as those found in Seidel~\cite{Se}, together with an adaptation of the ECH gluing theorem from Hutchings-Taubes~\cite{HT1,HT2}.

\subsection{ECH-type indices}
\label{subsection: ECH-type indices}

A key technical ingredient is the definition of ECH-type indices which is carried out in \cite{CGH2}.  We have already explained that the ECH differential counts ECH index $I_{ECH}=1$ holomorphic curves in the symplectization $\R\times N$.  We can also define an ECH-type index $I_{HF}$ in the Heegaard Floer situation so that the holomorphic curve count in the definition of the Heegaard Floer differential is precisely a count of $I_{HF}=1$ holomorphic curves in $\R\times [0,1]\times S$.  As in the ECH case, $I_{HF}$ is roughly the sum of two terms: the Fredholm index and the self-intersection number.  Moreover, we can define an ECH index $I_{W_+}$ for curves on $W_+$ (and $I_{\overline{W}_-}$ for curves on $\overline{W}_-$, defined later).

\section{The map $\Psi: PFC_{2g}(N) \to \widehat{CF}(S,\mathbf{a},h(\mathbf{a}))$}

In this section and the next, we explain the proof of the following:

\begin{thm}\label{thm: main}
The map $\Phi_* : \widehat{HF}(S,\mathbf{a},h(\mathbf{a})) \to PFH_{2g} (N)$ is an isomorphism which takes the Heegaard Floer contact invariant $c(\xi_{(S,h)})$ to the ECH contact invariant for $\xi_{(S,h)}$.
\end{thm}

\begin{rmk}
The fact that $\Phi_*$ is an isomorphism (i.e., that we do not need to consider orbit sets which intersect a page more than $2g$ times) is consistent with the adjunction inequality in knot Heegaard Floer homology, which states that if $S$ is a Seifert surface for a knot $K\subset M$, then $\widehat{HFK}(M,K,\mathfrak{s})\not=0$ only if $|\langle c_1(\mathfrak{s}),S\rangle| \leq 2g$; see \cite[Theorem~5.1]{Osz3}.
\end{rmk}

Theorem~\ref{thm: main} is proved by exhibiting a chain map
$$\Psi:PFC_{2g}(N) \to \widehat{CF}(S,\mathbf{a},h(\mathbf{a})),$$
constructed from a cobordism $\overline{W}_-$ which is similar to, but more complicated than, $W_+$, and proving that the induced map $\Psi_*$ on homology is the inverse of $\Phi_*$.

\subsection{A complication}

Let $B_-=(\R\times (\R/2\Z))- B_-^c$, where $B_-^c$ is $(-\infty,-2]\times[1,2]$ with the corners rounded; see the right-hand side of Figure~\ref{figure: bases}.  The naive candidate for the cobordism for $\Psi$ would be $W_-=\pi^{-1}(B_-)$, where $\pi: \R \times N \to \R\times S^1$ is as before.  It turns out that this naive candidate for $\Psi$ does not work for index reasons.  If we stack $W_+$ on top of $W_-$ in order to apply the usual chain homotopy argument (as in the next section), then we find that the Fredholm indices of the relevant degree $2g$ multisections $u$ from $\vec{\bf y}$ to itself are $-2g$, instead of $0$, as we would like.  In order to correct the index, we take a copy of the fiber, apply a multiple connect sum with $u$, and require that the curve pass through a marked point ${\frak m}$.

\subsection{The symplectic cobordism $(\overline{W}_-, \overline{\Omega}_- )$}

The above discussion motivates the definition of $(\overline{W}_-,\overline{\Omega}_-)$ and the point constraint ${\frak m}$. (The symplectic cobordism $(\overline{W}_+,\overline{\Omega}_+)$ can be defined similarly.)

Let $\overline{S}$ be the closed surface obtained by attaching a disk $D=\{\rho\leq 1\}$ with polar coordinates $(\rho,\phi)$ to $S$ along $\partial S$.  We extend $h:S\stackrel\sim\to S$ to $\overline{h}:\overline{S} \stackrel\sim\to \overline{S}$ so that all closed orbits of the suspension $\pi:\overline{N}\to S^1$ of $(\overline{S},\overline{h})$, contained in the open solid torus $int(D)\times S^1$, have intersection number $m\gg 2g$ with a fiber, with the exception of the orbit $\delta_0=\{\rho=0\}$. Note that $\overline{N}$ is obtained from $M$ by $0$-surgery along the binding. We extend the area form $\omega$ on $S$ to an area form $\overline{\omega}$ on $\overline{S}$, the stable Hamiltonian structure $(dt,\omega,\bdry_t)$ on $N$ to the stable Hamiltonian structure $(dt,\overline{\omega},\bdry_t)$ on $\overline{N}$, and the arcs $a_i$, $1\leq i\leq 2g$, to arcs $\overline{a}_i$ so that they start and end at the center $z_\infty$ of $D$ and restrict to radial lines in $D$.

We then define $\overline{W}_- =\pi^{-1} (B_-)$, where $\pi : \R \times \overline{N} \to \R \times S^1$. The symplectic form $\overline{\Omega}_-$ is the restriction of $ds\wedge dt +\overline{\omega}$ to $\overline{W}_-$ and the Lagrangian submanifolds $L_{{\bf \overline{a}}}$ are obtained by parallel transporting ${\bf \overline{a}}$ along $\partial B_-$ using the symplectic connection.

Finally, we pick the marked point ${\frak m}\in \overline{W}_-$ which is close to $\{z_\infty\}\times
B_-$ and in a generic position.

\subsection{The chain map $\Psi$}

Let $\overline{J}_-$ be an almost complex structure on $(\overline{W}_-,\overline{\Omega}_-)$ which is the restriction of an adapted almost complex structure on $\R\times\overline{N}$. Let $\gamma$ be an orbit set in $PFC_{2g}(N)$ and let $\vec{\bf y}$ be a $2g$-tuple of chords in $\widehat{CF}(S,{\bf a}, h({\bf a}))$.  If we write $\Psi(\gamma)= \sum_{\vec{\bf y}}  \langle \Psi (\gamma), \vec{\bf y} \rangle \vec{\bf y}$, then $\langle \Psi (\gamma), \vec{\bf y} \rangle$ counts degree $2g$, ECH index $I_{\overline{W}_-}=2$, $\overline{J}_-$-holomorphic multisections $u$ of the fibration $\pi: \overline{W}_-\to B_-$ which satisfy the following:
\begin{itemize}
\item $u$ is positively asymptotic to $\gamma$ and negatively asymptotic to $\vec{\bf y}$;
\item $u$ passes through the marked point ${\frak m}$ and has algebraic intersection number one with $\{z_\infty\}\times B_-$;
\item the boundary of $u$ is contained in distinct components of the Lagrangian boundary $L_{\bf a}$.
\end{itemize}
We will refer to such curves $u$ as {\em $\overline{W}_-$-curves} of ECH index $I_{\overline{W}_-}=2$.  Here the ECH index $I_{\overline{W}_-}$ is the index before taking into account the point constraint ${\frak m}$.

\begin{thm} \label{thm: Psi chain map}
The map $\Psi$ is a chain map.
\end{thm}

The proof of Theorem~\ref{thm: Psi chain map} requires more work than the proof of Theorem~\ref{thm: Phi chain map}, due to two factors: (i) $L_{\overline{\bf a}}$ is singular along $\{z_\infty\}\times \bdry B_-$ and (ii) the $\overline{W}_-$-curves are not contained in $W_-$.  Hence the limit of a sequence of $\overline{W}_-$-curves of ECH index $I_{\overline{W}_-}=3$ which pass through ${\frak m}$ can be more complicated; nevertheless, a more careful analysis of the limit curves yields Theorem~\ref{thm: Psi chain map}.

\section{Completion of the proof of Theorem~\ref{thm: isomorphism}}

We sketch the remaining steps in the proof of Theorem~\ref{thm: isomorphism}.

\subsection{Chain homotopy}
\label{subsection: chain homotopy}

In order to establish Theorem~\ref{thm: main}, it remains to prove that:

\begin{thm}\label{thm: identity}
The compositions $\Phi \circ \Psi$ and $\Psi \circ \Phi$ are chain homotopic to the identity.
\end{thm}

We consider the gluings of the two cobordisms $\overline{W}_+$ and $\overline{W}_-$ in different orders. In the case of $\Phi \circ \Psi$, we obtain a symplectic fibration over a cylinder with an open disk removed; in the case of $\Psi \circ \Phi$, we obtain a symplectic fibration over an annulus with one puncture (or, equivalently, with a strip-like end) on each boundary component.

\begin{figure}[ht]
\begin{overpic}[width=8cm]{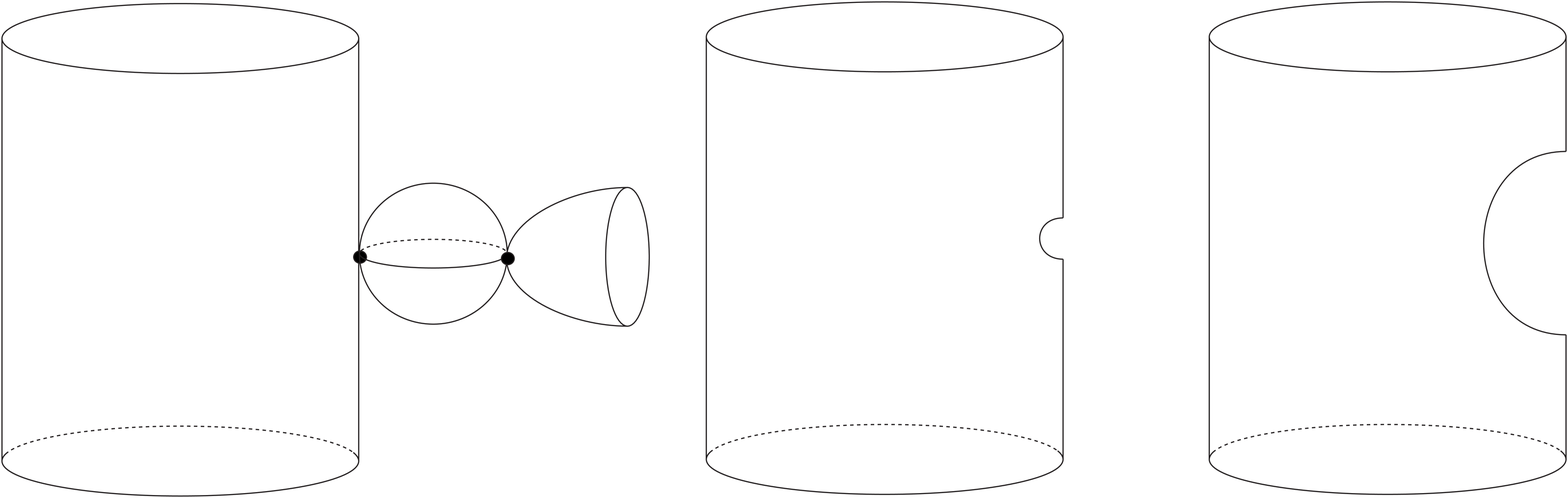}
\put(26.5,19.25){\tiny $\times$} \put(66.6,21){\tiny $\times$} \put(98.6,25){\tiny $\times$}
\end{overpic}
\caption{Degeneration (from right to left) of the cobordism for $\Phi \circ \Psi$. The location of (the projection of) the marked point is indicated by $\times$.}
\s
\label{figure: base degeneration 2}
\end{figure}

We sketch the proof for the composition $\Phi \circ \Psi$.  We degenerate the base to a nodal surface:
$$B_{-\infty}=(\R\times S^1)\sqcup S^2 \sqcup D^2/\sim,$$
where $p\in \R\times S^1$, $q_1\not=q_2 \in S^2$ and $r\in D^2$, and we identify $p\sim q_1$, $q_2\sim r$. (See the left-hand side of Figure~\ref{figure: base degeneration 2}.) The cobordisms corresponding to $\R\times S^1$, $S^2$ and $D^2$ are $\R\times \overline{N}$, $\overline{S}\times S^2$ and $\overline{S}\times D^2$; they are glued to give the fibration $\pi: \overline{W}_{-\infty}\to B_{-\infty}$. The marked point ${\frak m}$ is placed on the section  $\{z_\infty\}\times S^2$ of $\overline{S}\times S^2$. The holomorphic curves that we are counting then degenerate to the gluing of the following three types of curves:
\begin{enumerate}
\item trivial cylinders $u_1$ over closed orbits in $\R\times \overline{N}$;
\item curves $u_2$ in the class $[\overline{S} ] +2g[S^2 ] \in H_2 (\overline{S} \times S^2;\Z)$ which pass through (i) ${\frak m}$, (ii) $2g$ points in $\pi^{-1}(q_1)$ and (iii) each of the arcs $\overline{a}_i \times \{ q_2\}$; and
\item constant sections $u_3$ of $\overline{S}\times D^2$ with Lagrangian boundary $\overline{\bf a}\times \bdry D^2$.
\end{enumerate}
Here (ii) and (iii) come from the gluing constraints: (ii) from gluing $u_2$ to trivial cylinders $u_1$  and (iii) from gluing $u_2$ to constant sections $u_3$ with Lagrangian boundary $\overline{\bf a}\times \bdry D^2$. By a Gromov-Witten type computation, the count of (2) is $1$, implying the chain homotopy of $\Phi\circ\Psi$ with $id$.  A similar degeneration can be constructed for $\Psi \circ \Phi$.

\subsection{Stabilization}

Now that we have the isomorphism $$(\Phi_{(S,h)})_*:\widehat{HF}(S,{\bf a},h({\bf a}))\stackrel\sim\rightarrow PFH_{2g}(N),$$ where $\Phi_{(S,h)}$ refers to the $\Phi$-map corresponding to $(S,h)$, it remains to show that the map $PFH_{2g}(N)\to \lim \limits_{\longrightarrow \atop i}PFH(N)$ is an isomorphism.  This is proved in \cite{CGH2} by applying two positive stabilizations (corresponding to the connected sum with a trefoil knot) to $(S,h)$ to obtain $(S',h')$. We then complete ${\bf a}$ to ${\bf a'}$ by adding two extra basis arcs.  A comparison of $(\Phi_{(S,h)})_*$ and
$$(\Phi_{(S',h')})_*: \widehat{HF}(S',{\bf a'},h'({\bf a'}))\stackrel\sim\rightarrow PFH_{2g+2}(N)$$
implies that the map $(\mathfrak{I}_{2g})_*: PFH_{2g}(N)\rightarrow PFH_{2g+2}(N)$ is an isomorphism.  We similarly prove that $(\mathfrak{I}_i)_*$ is an isomorphism for $i\geq 2g$.

\s\n {\em Acknowledgements.}
Vincent Colin was supported by the Institut Universitaire de France, ANR Symplexe, and ANR Floer Power. Paolo Ghiggini was supported by ANR Floer Power. Ko Honda was supported by NSF Grant DMS-0805352. We are indebted to Michael Hutchings for many helpful conversations and for our previous collaboration which was a catalyst for the present work. We also thank Tobias Ekholm, Dusa McDuff, Ivan Smith and Jean-Yves Welschinger for illuminating exchanges.  Part of this work was done while KH and PG visited MSRI during the academic year 2009--2010. We are extremely grateful to MSRI and the organizers of the ``Symplectic and Contact Geometry and Topology'' and the ``Homology Theories of Knots and Links'' programs for their hospitality; this work probably would never have seen the light of day without the large amount of free time which was made possible by the visit to MSRI.

\end{document}